\documentclass[epsfig,12pt]{article}
\usepackage{amsfonts}
\usepackage{epsfig}
\usepackage{color}
\usepackage{graphicx}
\def\R{\mathbb{R}}

\newtheorem{definition}{{\sc Definition}}
\newtheorem{proposition}{{\sc Proposition}}

\newtheorem{remark}{Remark}

\date{November 17, 2012}

\begin{document}

\pagestyle{plain}
\rm

\title{ \Large {\bf To Split or Not to Split, That Is the Question in Some Shallow Water Equations}}

\author{
Vicente Mart\'{\i}nez\thanks{
    {\it E-mail address:}
                     martinez@mat.uji.es.}
    \\ {\small \it Departamento de Matem\'aticas \& Instituto de Matem\'aticas y sus Aplicaciones de Castell\'on (IMAC)} 
   \\ {\small \it  Universitat Jaume I, Campus de Riu Sec, 12071 Castell\'o, Spain}}

\maketitle

\begin{abstract}
 In this paper we analyze the use of time splitting techniques for solving shallow water equation. We discuss some
  properties that these schemes should satisfy so that interactions between the source term and the shock waves are controlled.   This paper shows
  that these schemes must be well balanced in the meaning expressed by   
  Greenberg and Leroux \cite{gre}. More specifically,  we analyze in what cases it is enough to verify an
{\it Approximate C-property} and in which cases it is required to verify an {\it Exact C-property} (see \cite{ber1},   \cite{ber2}).  
  We also include some numerical tests in order to justify our reasoning.

 \end{abstract}



\section{\bf  Introduction}
\hspace{0.2in} 
In this paper, our interest is to analyze time splitting schemes on conservation laws with source terms, also called
balance laws.  A prototype,  in one space dimension and under certain regularity hypotheses, is given by  the following system of partial differential equations
\begin{eqnarray} \label{ec1}
    \left\{ \begin{array}{l}
     W(x,t)_t + F(W(x,t))_x = G(x,W(x,t)) \:, \:\: (x, t) \in \R  \times \R^+ ,\\[0.25 cm]
     W(x,0)=W_0(x) \:, \:\: x \in \R ,  \end{array} \right.
      \end{eqnarray}
where $\: W: \mathbb{R} \times \mathbb{R}^+ \to \mathbb{R}^m \:$
is the vector of conserved variables,  $\: F: \mathbb{R}^m
\to \mathbb{R}^m  \:$ is the vector of fluxes and $\:G:  \mathbb{R}^{m+1}
\to \mathbb{R}^m  \:$ is the source term.

\par
Recently, there has been some controversy related to the application of time splitting techniques on hyperbolic equations involving solutions with shock waves, as is the case of shallow water equations, which are introduced later. 
Usually,  a time splitting numerical schemes to solve  (\ref{ec1}) consists of solving consecutively the homogeneous equation
  \begin{eqnarray} \label{ec2}
     W(x,t)_t + F(W(x,t))_x = 0,
      \end{eqnarray}

\noindent
and  the ordinary differential equation 
\begin{eqnarray} \label{ec3}
     W(x,t)_t = G(x,W(x,t)).
      \end{eqnarray}

\par
LeVeque notices in \cite{lev1} that such schemes can easily fail by the presence of shock waves in solving
 (\ref{ec2}). These shock waves involve large changes in the solution which can not be captured in solving
  (\ref{ec3}).

\par
On the other hand, some authors such as  Ma, Sun and Yin (see \cite{ma}) use a time integrating scheme with two-step predictor-corrector sequence quite successfully.
Striba use splitting techniques (see \cite{ski}) in meteorology  models on the term that represent de Coriolis acceleration.
In addition, Wicker and Skamarock (see \cite{wic}) use time-splitting methods for integrating the elastic equations.

\par
In 1994, Berm\'udez and V\'azquez  (see \cite{ber1},  \cite{ber2}) introduce the concept of 
{\it Exact C-property} and 
{\it Approximate C-property} in order to identify numerical schemes with an acceptable level of accuracy in the resolution of shallow water equations (well-balanced scheme). Going deeper into 
well-balanced scheme idea, we can find the work of
Greenberg and Leroux \cite{gre}, in which they propose a numerical scheme that preserves a balance 
 between the source terms and internal forces due to the presence of shock waves.

\par
Another point of view, more recently, it is provided by Lubich (see \cite{lub}), who gives an error analysis of Strang-type splitting integrators for nonlinear Schr\"{o}dinger equations. Holdahl, Holden and Lie use an adaptive grid refinement and 
a shock tracking technique to construct a front-tracking method for hyperbolic conservation laws. They combine the operator splitting to study the shallow water equations (see \cite{holda}). Holden, Karlsen, Risebro and Tao
 (see\cite{holden2}) show that the Godunov and Strang splitting methods converge with the expected rates if the initial data are sufficiently regular. Finally, in this way the reader can find a deep study of splitting methods for partial differential 
 equations in \cite{holden1}, where some analysis of conservation and balance laws are included.

\par
 We base our analysis on the ideas presented in 
  \cite{ber1},  \cite{ber2} and  \cite{gre}. 
 From these studies, it follows that the numerical scheme used in solving (\ref{ec2}) and the numerical scheme used in solving (\ref{ec3}) cannot be whatever, even though they have a high degree of accuracy. These must be balanced so that interactions between the source term in (\ref{ec3}) and the shock waves in (\ref{ec2}) are controlled. 
In this framework, we will analyze conditions to be verified by splitting schemes in order to avoid spurious oscillations, which are created in this type of equations. More specifically, in which cases it is enough to verify an
{\it Approximate C-property} and in which cases it is required to verify an {\it Exact C-property}.

\par
The rest of this paper is organized as follows: in Section 2, we introduce the 
governing equations for the one dimensional shallow water model. We also analyze two kinds of time splitting schemes
to identify which conditions must satisfy a scheme for solving this type of equations. In Section 3, we show the four test 
problems  that we will use to test the performance of the schemes described above. In Section 4, we present the
     numerical results we have obtained. Finally, in Section 5 we reason the final conclusions we have reached by using
        the results obtained with the presented schemes.

   \vspace{0.2cm}
\section{\bf The one dimensional shallow water equations}

\subsection{\bf Governing equations}

 \hspace{0.2in}
In this section we consider Eq. (\ref{ec1}) with

\begin{equation}\label{ec4}
\begin{array}{c}
    W=\left(%
    \begin{array}{c}
      h \\
      q \\
   \end{array}%
    \right),
    \hspace{1cm} 
    F(W)=\left(%
      \begin{array}{c}
        q \\[0.4 cm]
     \displaystyle   \frac{q^2}{h}+\frac{1}{2}gh^2 \\
      \end{array}%
    \right), \hspace{1cm} 
G(x,W)= \left(%
\begin{array}{c}
  0 \\
  -ghb'(x)    \\
    \end{array}
   \right),
\end{array}
\end{equation}

\noindent
where the unknowns of the problem are: the water height is $ h $ and the
flow per unit length is $ q = hu .$  Here $ u $ is the average
vertical speed in the direction of the axis $ x $ (see Fig.1).  $ F $ is the flux of 
conservative variables and $ g = 9.81ms^{-2} $ is the acceleration
due to gravity. The source term $ G $
 models the bottom variation given by the function $ b (x). $

%%%% figura 1%%%%%%%%%%%%%%%%%%%%%%%%%%%%%%%%%%%%%%%
% \vspace{2 cm}
  \begin{figure}[ht]
   \begin{center}
    \begin{picture}(500,180)
    
     \put(70,-40){ \includegraphics[width=10cm]{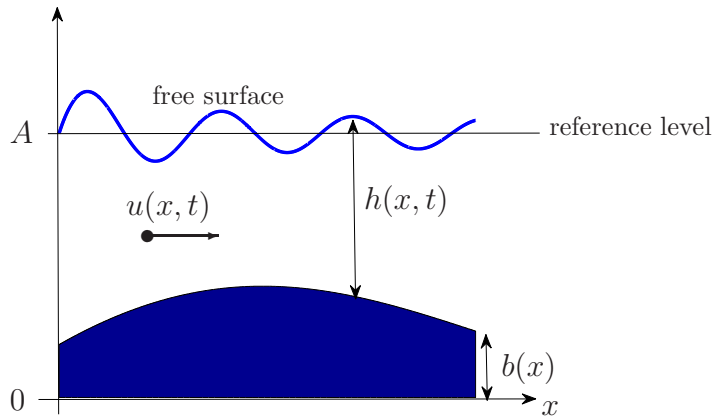} }     
    \put(170,120){\footnotesize free surface}     
       \put(320,108){\footnotesize reference level}  
           
     \put(170,70){\vector (1,0){25}}
       \put(160,78){$u(x,t)$}
        \put(165,67){$\bullet$}
     
       \put(116,4){$0$}
       \put(116,105){$A$} 
        \put(302,17){$b(x)$}
      \put(250,80){$h(x,t)$}
      \put(317,2){$x$}
        \end{picture}
    \end{center}
    \caption{Shallow water variables.}
    \end{figure}
%%%%%%%%%%%%%%%%%%%%%%%%%%%%%%%%%%%%%%%%%%%%%    

\vspace{0.5 cm} 
Let us consider numerical solvers based upon the decomposition $F(W)_x=A(W) \: W_x,$ where
\begin{equation}\label{ec5}
\begin{array}{c}
    A(W)= \left(%
    \begin{array}{cc}
      0 & 1 \\[0.4 cm]
        \displaystyle   - \frac{q^2}{h^2}+ gh & \displaystyle  2 \frac{q}{h}  \\
   \end{array}%
    \right)  \: \mbox{is the Jacobian matrix of } \: F(W).
   \end{array}
\end{equation}

\vspace{1cm}

So that, system (\ref{ec1})-(\ref{ec4}) is hyperbolic $(h>0),$ then  $A=X \Lambda X^{-1},$ where

\begin{equation}\label{ec6}
\begin{array}{c}
    \Lambda =\left(%
    \begin{array}{cc}
      \lambda_1 & 0 \\
       0 & \lambda_2 \\
   \end{array}%
    \right),
    \hspace{1cm} 
    X=\left(%
      \begin{array}{cc}
     1  & 1 \\
       \lambda_1 & \lambda_2 \\
   \end{array}
    \right), \hspace{1cm} 
X^{-1}= \frac{1}{\lambda_2 - \lambda_1} 
   \left(
      \begin{array}{rr}
     \lambda_2  & -1 \\
   - \lambda_1 & 1 \\
   \end{array}%
    \right), 
\end{array}
\end{equation}

\noindent
where  $\displaystyle \lambda_1= \frac{q}{h}+ \sqrt{gh}$ and    
$\displaystyle \lambda_2= \frac{q}{h} - \sqrt{gh}.$

\vspace{0.75 cm}

Our analysis needs the definition of some conservation properties  given by Berm\'udez and V\'azquez in \cite{ber2}.
Since only the source term involves the bed slope, an inadequate choice of the numerical schemes can give
spurious oscillations. So, these conservation properties try to identify which kind of schemes have good 
behavior in equations with source term. Berm\'udez and V\'azquez characterize the good behavior of the numerical scheme in the manner in which the scheme approximates a steady solution representing the state of water at rest.
They introduce the stationary problem ({\it Problem SP}) given by 
 $q(x, t)=0$ and $h(x, t)=H(x)$ and define the following conservation properties:

\begin{definition}  {\bf Exact C-property}.
We say that a scheme satisfies the Exact C-Property if it is exact when applied to the Problem SP.
\end{definition}

\begin{definition}  {\bf Approximate C-property}.
We say that a scheme satisfies the Approximate C-Property if it is accurate to the order 
$\Theta (\Delta x^2)$ when applied to Problem SP.
\end{definition}

When a numerical scheme does not satisfy any of these conservation properties then the propagation of spurious
 oscillations is also present in non stationary problems.

\subsection{\bf Central numerical schemes}

The differential formulation of the homogeneous equation given in (\ref{ec2}) does not  admit discontinuous solutions. 
 So, since these solution are physically relevant in this context, we need a suitable formulation of the problem to support discontinuous solutions. In this sense, it is usual to consider the following integral formulation

\begin{equation}\label{integ}
\int(W dx-F(W) dt)=0.
\end{equation}

\par
The numerical schemes use usually (\ref{integ}) in order to approximate (\ref{ec2}). To do that, it is introduced 
 a control volume
 in the space $(x, t)$ of
dimensions $\Delta x \times \Delta t$. Next, it is  evaluated the integral
(\ref{integ}) in this volume control
\begin{displaymath}
\int_{x_{j-1/2}}^{x_{j+1/2}}\big(W(x,t^{n+1})-W(x,t^n)\big)dx+
\int_{t^n}^{t^{n+1}}\big(F(W(x_{j+1/2},t))-F(W(x_{j-1/2},t))\big)dt=0.
\end{displaymath}

\par
\par

Dividing by $\Delta x$ we obtain

\begin{eqnarray*}
\frac{1}{\Delta x}
\int_{x_{j-1/2}}^{x_{j+1/2}}W(x,t^{n+1})dx&=&\frac{1}{\Delta x}
\int_{x_{j-1/2}}^{x_{j+1/2}}W(x,t^n)dx\\&-&\frac{\Delta t}{\Delta
x}\Big[\frac{1}{\Delta t}\int_{t^n}^{t^{n+1}}F(W(x_{j+1/2},t))dt
-\frac{1}{\Delta t}\int_{t^n}^{t^{n+1}}F(W(x_{j-1/2},t))dt\Big].
\end{eqnarray*}

\par

Thus, we deduce the conservation formula 
\begin{equation}\label{forcon1}
\overline{W}_j^{n+1}=\overline{W}_j^n-\frac{\Delta t}{\Delta
x}[F_{j+1/2}-F_{j-1/2}],
\end{equation}

\noindent
where $\overline{W}_j^n$ is an average
\begin{displaymath}
\overline{W}_j^n=\frac{1}{\Delta
x}\int_{x_{j-1/2}}^{x_{j+1/2}}W(x,t^n)dx
\end{displaymath}
at time $t=t^n$ inside the interval
\begin{displaymath}
I_j=[x_{j-1/2},x_{j+1/2}]
\end{displaymath}
whose length is
\begin{displaymath}
\Delta x=x_{j+1/2}-x_{j-1/2}.
\end{displaymath}

\par
The flux in (\ref{forcon1}) can be interpreted as the average in time
of the physical flux, i.e.,
\begin{equation}\label{flujonum}
F_{j+1/2}\approx\frac{1}{\Delta t}\int_{t^n}^{t^{n+1}}F(W(x_{j+1/2},t))dt.
\end{equation}

\par
 Conservative numerical methods for (\ref{ec2}) are based
in (\ref{forcon1}), and they are determined by the expression of
the numerical flux $F_{j+1/2}$.

\par
The time splitting numerical schemes (see \cite{tor3}), in each time step, act as follows: taking into account  the 
    initial condition  $W^{n},$ we  solve (\ref{ec2})  and we obtain  $\widehat{W}^{n+1}.$ Then, by using the 
    initial condition     $\widehat{W}^{n+1},$ we solve    (\ref{ec3})  and  obtain $W^{n+1}.$

\vspace{0.25 cm}

More specifically, if we use the operators

  $\bullet$  ${\cal A} (W^{n})=W^{n+1}  \: \:  \mbox{such as }$ 
    $$W^{n} \rightarrow 
    \left\{ \begin{array}{l}
     W_t + F(W)_x =0  , \: (x, t) \in \R  \times [t^n, t^{n+1}] ,\\[0.25 cm]
     W(x,t^n)=W^n ,  \: x \in \R.  \end{array} \right\}  \rightarrow \widehat{W}^{n+1} .$$

 $\bullet$  ${\cal S} (W^{n})=W^{n+1}  \: \:  \mbox{such as }$ 
 $$W^{n} \rightarrow 
    \left\{ \begin{array}{l}
     W_t  = G(x,W) , \: (x, t) \in \R  \times [t^n, t^{n+1}] ,\\[0.25 cm]
     W(x,t^n)= \widehat{W}^{n+1}  ,  \: x \in \R.    \end{array} \right\}  \rightarrow W^{n+1} .$$

\par

Thus,

      $${\cal S}^{(\Delta t)} \cdot  {\cal A}^{(\Delta t)}  (W^{n})=W^{n+1} .$$

\par
In the first step of the splitting procedure, we solve  (\ref{ec2}) at each time step. To do that and for the sake of simplicity, we consider the $Q$-scheme of van Leer (see \cite{vaz}) which uses a matrix $Q$ satisfying some properties and the numerical fluxes
  \begin{equation}\label{ec10}
F_{j-1/2} = \phi(W_{j-1}^n,W_j^n) , \: \: \:  \: \: \:   F_{j+1/2} = \phi(W_{j}^n,W_{j+1}^n),
\end{equation} 

\noindent
where the numerical flux function $\phi$ is given by
 \begin{equation}\label{ec11}
 \phi(U,V)=\frac{F(U)+F(V)} {2} - \frac{1}{2} |Q(U,V)| (V-U).
  \end{equation}  
  
  \par
  
  A possible choice of the matrix $Q$ can be the Jacobian (\ref{ec5}) of the system (\ref{ec4}) evaluated at the arithmetic
   mean, i.e.
  \begin{equation}\label{ec12}
 Q(U,V)=A\left( \frac{U+V} {2} \right).
  \end{equation}    
  
  So that, we have
   \begin{equation}\label{ec13}
              \left| Q(W_{j\pm\frac{1}{2}}^n) \right| =
                            X(W_{j\pm1}^n, W^n_{j}) \left| \Lambda(W_{j\pm1}^n, W^n_{j}) \right|
                                X^{-1}(W_{j\pm1}^n, W^n_{j}) ,
  \end{equation}    
  
  \noindent  
  where   $X(W_{j\pm1}^n, W^n_{j}), \left| \Lambda(W_{j\pm1}^n, W^n_{j}) \right|
        \:\: \mbox{and} \:\:
                                X^{-1}(W_{j\pm1}^n, W^n_{j}) $ 
    are  evaluated at  $\: \displaystyle \frac{W_{j\pm1}^n + W^n_{j}}{2}.$

   \vspace{0.4 cm}
          
  \par                             
So, we obtain the numerical fluxes
   \begin{equation}\label{ec14}
F_{j \pm 1/2} = \frac{F(W_{j\pm1}^n)+F(W^n_{j})} {2} - 
                        \frac{1}{2} \: \left| Q(W_{j\pm\frac{1}{2}}^n) \right| \: (\mp W^n_{j} \pm W_{j\pm1}^n).
\end{equation}

  \par  
  The second step of the procedure is to solve  (\ref{ec3}). To do that, we get the solutions of the homogeneous
   equation 
   (\ref{ec2}):   $ \hat{h}(x_j,t^{n+1})$ and  $\hat{q}(x_j, t^{n+1})$ 
and we  solve the following initial value ODE problem for each $x_j$

 \begin{eqnarray} \label{ec14_2}
     \left\{ \begin{array}{l}
  \displaystyle \frac{d}{dt} \left( \begin{array}{l}
     h(x_j, t) \\[0.25 cm]
     q(x_j, t) \end{array} \right) =   \left( \begin{array}{c}
      0 \\[0.15 cm]
     - g h(x_j,t)  b'(x_j) \end{array} \right); \: \: \: \: \:  t\in [t^n, t^{n+1}] , \\[0.65 cm]
      \left( \begin{array}{l}
     h(x_j,t^n) \\[0.15 cm]
     q(x_j, t^n) \end{array} \right) =       
      \left( \begin{array}{l}
     \hat{h}(x_j,t^{n+1}) \\[0.25 cm]
     \hat{q}(x_j, t^{n+1}) \end{array} \right)      .
     \end{array} \right.
      \end{eqnarray}

  \par
  
   The first equation has the solution:    

\begin{equation}\label{ec15}
h(x_j,t^{n+1})=   \hat{h}(x_j,t^{n+1}).
\end{equation}

 \par  
 
 To solve the second one,  we calculate
    
     $$ \displaystyle \frac{d q(x_j, t) }{dt}  =    - g h(x_j, t)  b'(x_j),$$
       
          $$  q(x_j, t)   =    - \int_{t^n}^{t} g h(x_j, s)  b'(x_j) \: ds ; \:\:   t\in   (t^n, t^{n+1}].  $$
   
         \vspace{0.2 cm }    
Therefore,
        \vspace{-0.2 cm }    
     $$  q(x_j, t^{n+1})   =   \hat{q} (x_j, t^{n+1})   -   g \cdot  b'(x_j) \cdot     
              \int_{t^n}^{t^{n+1}}  h(x_j, s)  \: ds .  $$
         
         \par
Finally,  by using a trapezoidal rule, we obtain 

\begin{equation}\label{ec16}
  q(x_j, t^{n+1})   =   \hat{q} (x_j, t^{n+1})   -   g \cdot  b'(x_j) \cdot     
              \frac{\Delta t}{2} (  h(x_j, t^n) +    \hat{h}(x_j, t^{n+1}) ) .  
 \end{equation}
      
   \begin{remark} 
   We denote the time splitting scheme given by 
   (\ref{ec14})-(\ref{ec15})-(\ref{ec16}) as {\bf Q-tra1}.
  \end{remark}
                                   
 \begin{proposition} 
   The numerical scheme {\bf Q-tra1} satisfies an approximate C-property.
  \end{proposition} 

{\bf Proof}

Problem SP satisfies 
$\lambda_1=\sqrt{gh} \: $ and $ \: \lambda_2= - \sqrt{gh}.\:$ In addition (see Fig. 1), $b(x)+h(x,t)=b(x)+H(x)=A.$
Therefore, 
  \begin{equation}\label{ec100} b(x_j)=A-H(x_j).
   \end{equation}

      In order to solve  (\ref{ec2}) with (\ref{ec14}), we compute
      
   $$ \left| Q(W_{j + \frac{1}{2}}^n) \right| = \frac{1}{\lambda_2 - \lambda_1}
           \left( \begin{array}{cc}
               |\lambda_1| \lambda_2 - |\lambda_2 | \lambda_1 & -|\lambda_1| +  |\lambda_2|  \\[0.25 cm]
                \lambda_1 \lambda_2 ( |\lambda_1| - |\lambda_2| ) &  - \lambda_1 | \lambda_1| + \lambda_2  |\lambda_2|
             \end{array} \right)  = 
            \left( \begin{array}{cc} \displaystyle
              \sqrt{g \: \frac{h_{j+1} + h_j}{2}}     &  0  \\[0.45 cm]
               0           &   \displaystyle \sqrt{g \: \frac{h_{j+1} + h_j}{2}}  \end{array} \right)   .           $$
      
     Then,  
      
      $$F^n_{j + \frac{1}{2}}= \frac{1}{2}
            \left( \begin{array}{c} \displaystyle
             (h_{j+1} - h_j) \: \sqrt{g \: \frac{h_{j+1} + h_j}{2}}   \\[0.45 cm]
           \displaystyle   \frac{1}{2} \: (h^2_{j+1} + h^2_j) \end{array} \right)   .           $$

       Analogously, we obtain             
                               $$F^n_{j - \frac{1}{2}}= \frac{1}{2}
            \left( \begin{array}{c} \displaystyle
             (h_{j} - h_{j-1}) \: \sqrt{g \: \frac{h_{j} + h_{j-1}}{2}}  \\[0.45 cm]
           \displaystyle   \frac{1}{2} \: (h^2_{j} + h^2_{j-1}) \end{array} \right)   .           $$

   Finally, using (\ref{forcon1}), we have  the numerical approximation of Eq. (\ref{ec2})

$$\widehat{W}_j^{n+1} = \left( \begin{array}{c}
                                      \hat{h}_j \\[0.5 cm]
                                    \hat{q}_j \end{array} \right) = 
                                     \left( \begin{array}{c}
                                      {h}_j + \displaystyle  \frac{\Delta t}{2 \: \Delta x}
          \left(  (h_{j+1}-h_j) \sqrt{g  \: \frac{h_{j+1}+h_j}{2} }  - (h_{j}-h_{j-1}) \sqrt{g  \: \frac{h_{j}+h_{j-1}}{2} } \right)\\[0.5 cm]
           \displaystyle   -  \frac{g  \: \Delta t}{4  \: \Delta x}  \left( h_{j+1}^2-h_{j-1}^2    \right) 
                                        \end{array} \right)      $$
                               
     \par
                                
    Next, to have the numerical approximation of Eq. (\ref{ec3}),  we use  (\ref{ec15}) and (\ref{ec16}). By substituting, 
    we  obtain
        
        \begin{equation}\label{ec17}
                      h_j^{n+1}=   
   {h}_j + \displaystyle  \frac{\Delta t}{2 \: \Delta x}
          \left(  (h_{j+1}-h_j) \sqrt{g  \: \frac{h_{j+1}+h_j}{2} }  - (h_{j}-h_{j-1}) \sqrt{g  \: \frac{h_{j}+h_{j-1}}{2} } \right)
         \end{equation}
and    
          
          \begin{equation}\label{ec18}
               q_j^{n+1}=       \displaystyle   -  \frac{g  \: \Delta t}{4  \: \Delta x}  \left( h_{j+1}^2-h_{j-1}^2    \right)   -
                 g \cdot  b'(x_j) \cdot  \frac{\Delta t}{2} (  h(x_j, t^n) +    \hat{h}(x_j, t^{n+1}) ) .  
         \end{equation}

                      From  (\ref{ec17}), we deduce that  $h_j^{n+1}=   {h}_j $  is accurate to the order 
$\Theta (\Delta x^2).$   And from (\ref{ec18}) and (\ref{ec100}) we have that $ q_j^{n+1}= 0$ is exact when we take
$$ b'(x_j) = \frac{b(x_{j+1})- b(x_{j-1})}{2 \: \Delta x}. $$

Thus, the proof ends.
\begin{flushright}
$\square$
\end{flushright}

\vspace{0.75 cm}
  \begin{remark} 
   If we use another convergent quadrature rule to approximate the integral
  $\displaystyle \int_{t^n}^{t^{n+1}}  h(x_j, s)  \: ds ,$ then {\sc Proposition} 1 is also satisfied.
  \end{remark}
                              
   \vspace{0.75 cm}                            
                               
   \subsection{\bf Upwind numerical schemes}                            
                             
                  It is well-known  that conservation laws with source terms can be solved with high accuracy 
                  and discontinuities are well captured by using upwind schemes. For instance, you can see the works of
                   LeVeque and Yee \cite{lev2} or V\'azquez-Cend\'on \cite{vaz}.
                  Furthermore, this kind of schemes have bigger stability regions than schemes using centered 
                  approximations  of the source term (see \cite{ber2}). In this section, we will use some
                   approaches in order to build an upwind scheme using some ideas contained in \cite{vaz}
                   and making changes to solve (\ref{ec3}).

         \par 
                   In order to solve (\ref{ec3}), we need to introduce two numerical source functions, $G_L$ on the left and
                    $G_R$ on the right, to upwind the source term. We also use two matrices $D_L$ and $D_R$ 
                    with the aim of making clear the upwind process.                 
                       Below  we explain the process in detail.                    
                    
       \par
                     We will use, in each time step,
                   the average of solutions of  the following equations:                                           
                               
        \begin{equation}\label{ec19}
    \left\{ \begin{array}{l}
     W_t  = G_L(x,W) , \: (x, t) \in \R  \times [t^n, t^{n+1}] ,\\[0.25 cm]
     W(x,t^n)= \widehat{W}^{n+1}  ,  \: x \in \R.    \end{array} \right. 
        \end{equation}    
        
        and                           
               \begin{equation}\label{ec20}
    \left\{ \begin{array}{l}
     W_t  = G_R(x,W) , \: (x, t) \in \R  \times [t^n, t^{n+1}] ,\\[0.25 cm]
     W(x,t^n)= \widehat{W}^{n+1}  ,  \: x \in \R ,
        \end{array} \right. 
        \end{equation}                              
                               
                     where
               $$ G_L(x,W)=D_L( \overline{W}_L) \, G (x, W) \:\:     \mbox{and} \:\:
              G_R(x,W)=D_R( \overline{W}_R)  \, G (x, W);$$

     \par
   and functions $ \overline{W}_L$ and $ \overline{W}_R$ change in each computacional cell in the 
        discretization problem
     $$ \overline{W}_L=\frac{W_j+W_{j-1}}{2}  \:\:     \:\:   \mbox{and} \:\:   \:\:  
          \overline{W}_R=\frac{W_{j+1}+W_j}{2}. $$

         Thus, we obtain the solution of (\ref{ec19})-(\ref{ec20}), $ W_L^{n+1}$ and  $ W_R^{n+1}$        
         respectively. Then, in each time step, we take as solution of  (\ref{ec3}) the following average
         
                     $$W^{n+1}= \frac{W_L^{n+1}+W_R^{n+1}}{2}.$$

       In other words, for the discretization of the problem, we can substitute the ODE in  (\ref{ec14_2}) by

        $$    \displaystyle \frac{d}{dt} \left( \begin{array}{l}
                             h (x_j, t)\\[0.25 cm]
                             q (x_j,t)\end{array} \right) =   D_L(\overline{W}_{L_j}) \, G (x_j, W_j) =
                              \left( \begin{array}{cc}
                              d_{11}^L & d_{12}^L \\[0.15 cm]
                              d_{21}^L & d_{22}^L \end{array} \right) \cdot
                            \left( \begin{array}{c}
                                  0 \\[0.15 cm]
                                 - g h(x_j,t) b'(x_j) \end{array} \right) =
                           \left( \begin{array}{c}
      -  d_{12}^L g h_j b'_j \\[0.15 cm]
     - d_{22}^L g h_j b'_j \end{array}  \right).  $$

Therefore, instead of (\ref{ec14_2}), we have two
ODE         
                               
 \begin{eqnarray} \label{ec21}
     \left\{ \begin{array}{l}
  \displaystyle \frac{d}{dt} \left( \begin{array}{l}
     h_j \\[0.25 cm]
     q_j \end{array} \right) =    \left( \begin{array}{c}
      -  d_{12}^L g h_j b'_j \\[0.15 cm]
     - d_{22}^L g h_j b'_j \end{array}  \right); \: \: \: \: \:  t\in [t^n, t^{n+1}] , \\[0.65 cm]
      \left( \begin{array}{l}
     h(x_j,t^n) \\[0.15 cm]
     q(x_j, t^n) \end{array} \right) =       
      \left( \begin{array}{l}
     \hat{h}(x_j,t^{n+1}) \\[0.25 cm]
     \hat{q}(x_j, t^{n+1}) \end{array} \right)      .
     \end{array} \right.
      \end{eqnarray} 
      
       and  analogously      from (\ref{ec20}) we obtain        

    \begin{eqnarray} \label{ec22}
     \left\{ \begin{array}{l}
  \displaystyle \frac{d}{dt} \left( \begin{array}{l}
     h_j \\[0.25 cm]
     q_j \end{array} \right) =    \left( \begin{array}{c}
      -  d_{12}^R g h_j b'_j \\[0.15 cm]
     - d_{22}^R g h_j b'_j \end{array}  \right); \: \: \: \: \:  t\in [t^n, t^{n+1}] , \\[0.65 cm]
      \left( \begin{array}{l}
     h(x_j,t^n) \\[0.15 cm]
     q(x_j, t^n) \end{array} \right) =       
      \left( \begin{array}{l}
     \hat{h}(x_j,t^{n+1}) \\[0.25 cm]
     \hat{q}(x_j, t^{n+1}) \end{array} \right)      .
     \end{array} \right.
      \end{eqnarray} 
           
 where  $  d_{12}^R$ and $ d_{22}^R$ are coefficients corresponding to matrix $D_R.$

 Finally, we compute the solutions of (\ref{ec21}) and (\ref{ec22}), $W_{Lj}^{n+1}$ and $W_{Rj}^{n+1}.$  
                                       
 \begin{remark} 
   Matrices $D_L$ and $D_R$ must be chosen in coordination with the numerical method used in solving
  (\ref{ec2}) to get a well balanced scheme. For example, regarding  (\ref{ec14}), we can take  
   \begin{eqnarray} \label{ec23}
    D_L=  (I+|Q|Q^{-1})  \: \:   \mbox{and} \:\:  
            D_R=  (I-|Q|Q^{-1}).       
   \end{eqnarray}
   In this way we obtain consistency, since   $D_L+D_R=2I,$ we have
             $$ \frac{G_L + G_R}{2}= \frac{D_L G+ D_R G}{2}= \frac{(D_L+D_R) G}{2}=G.$$
   
  \end{remark}
                                   
                \begin{remark} 
   We denote the time splitting scheme given by 
   (\ref{ec14})-(\ref{ec21})-(\ref{ec22})-(\ref{ec23}) as {\bf Q-tra2}.
  \end{remark}

 \begin{proposition} 
   The numerical scheme {\bf Q-tra2} satisfies an exact C-property.
  \end{proposition} 

{\bf Proof}

Acting in a similar way to {\sc Proposition} 1, we have                               
                               
$$\widehat{W}_j^{n+1} = \left( \begin{array}{c}
                                      \hat{h}_j \\[0.5 cm]
                                    \hat{q}_j \end{array} \right) = 
                                     \left( \begin{array}{c}
                                      {h}_j + \displaystyle  \frac{\Delta t}{2 \: \Delta x}
          \left(  (h_{j+1}-h_j) \sqrt{g  \: \frac{h_{j+1}+h_j}{2} }  - (h_{j}-h_{j-1}) \sqrt{g  \: \frac{h_{j}+h_{j-1}}{2} } \right)\\[0.5 cm]
           \displaystyle   -  \frac{g  \: \Delta t}{4  \: \Delta x}  \left( h_{j+1}^2-h_{j-1}^2    \right) 
                                        \end{array} \right)   .   $$
                                        
                                        On the other hand, 
                                        $$D_1(\overline{W}_{L_j})= \left( \begin{array}{cc}
                                          1   & \frac{1}{\sqrt{g\frac{h_j + h_{j-1}}{2}}} \\
                                           \sqrt{g\frac{h_j + h_{j-1}}{2}}  & 1        
                                              \end{array} \right)   \: \: \mbox{and} \: \:
                                                  D_2(\overline{W}_{R_j})= \left( \begin{array}{cc}
                                          1   & - \frac{1}{\sqrt{g\frac{h_{j+1} + h_{j}}{2}}} \\
                                          -  \sqrt{g\frac{h_{j+1} + h_{j}}{2}}  & 1        
                                              \end{array} \right)      .         $$

                Now,  taking in   (\ref{ec21})           
                      $$ b'(x_j) = \frac{b(x_{j})- b(x_{j-1})}{ \Delta x} 
              \: \:    \mbox{ and}   \:\:  h(x_j, t)=    \frac{h_{j}+ h_{j-1}}{2} $$      to integrate, we have

  \begin{equation}\label{ec26}
                      h_{L_j}=   
   \hat{h}_j - \displaystyle  \frac{\Delta t}{ \Delta x}  (b_{j}-b_{j-1}) 
          \sqrt{g  \:  \frac{h_{j}+h_{j-1}}{2} }  
          \end{equation}
and    
          
          \begin{equation}\label{ec27}
                     q_{L_j}=   
   \hat{q}_j - \displaystyle  \frac{\Delta t}{ 2 \: \Delta x}  g (b_{j}-b_{j-1}) 
          (h_{j}+h_{j-1}) 
         \end{equation}

Acting in a similar way to  (\ref{ec22}), we have

 \begin{equation}\label{ec28}
                      h_{R_j}=   
   \hat{h}_j - \displaystyle  \frac{\Delta t}{ \Delta x}  (b_{j+1}-b_{j}) 
          \sqrt{g  \:  \frac{h_{j+1}+h_{j}}{2} }  
          \end{equation}
and    
          
          \begin{equation}\label{ec29}
                     q_{R_j}=   
   \hat{q}_j - \displaystyle  \frac{\Delta t}{ 2 \: \Delta x}  g (b_{j+1}-b_{j}) 
          (h_{j+1}+h_{j}) 
         \end{equation}

 In addition,   from (\ref{ec100}) we have $b_j=A-h_j,$ then 
  
 \begin{eqnarray*}
              h_j^{n+1}= \frac{h_{L_j}  +     h_{R_j}}{2}=   h_j 
         \:  \:  \:  \:  \:  \mbox{and}   \:  \:  \:  \:  \: 
       q_j^{n+1}= \frac{q_{L_j}  +     q_{R_j} }{2}=   0 .
     \end{eqnarray*}

This concludes the proof.
\begin{flushright}
$\square$
\end{flushright}

\section{Numerical experiments}
 
In this section, we will discuss four test problems for shallow water equations.
 The first test problem is the dam-break problem with variable topography, which 
 is one of the most basic problems with source term.
 The second one is a stationary problem, also with source term that represents a smooth bottom.
  The third one represents a tidal wave flow
propagating over an irregular topography. This test 
   was discussed by Berm\' udez and  V\'azquez \cite{ber2} and it is one of the most popular test to check the performance of a numerical scheme which tries to be effective in solving shallow water equations. It represents a severe test regarding to irregular topography and the long time over which is applied. Finally, the last test represents a numerical simulation of a tidal wave on the shoreline with friction effects, which introduces a wet/dry front.  
  
 \subsection*{Test 1: Dam-break with variable topography}

We consider the one dimensional shallow water equations (\ref{ec4}) on the domain 
  $(x, t) \in [0,1]  \times \R^+.$  The initial data are
   
  \begin{eqnarray*} 
  b(x)= 
      \left\{  \begin{array}{l}
   \displaystyle   \frac {1}{8} \cos\left(10 \pi \left(x-\frac{1}{2}\right) \right)+1 , \:  \frac{2}{5}<x<\frac{3}{5} , \\[0.45 cm]
    0, \mbox{otherwise;} \end{array} \right. 
        \:\:\:
     W(x, 0) =   \left\{ \begin{array}{l} 
           \left( \begin{array}{c}
     1- b(x) \\[0.15 cm]
     0  \end{array} \right) , \: x<\frac{1}{2} ,     \\[0.5 cm]
            \left( \begin{array}{c}
     0.5-b(x) \\[0.15 cm]
     0  \end{array} \right) ,  \: x>\frac{1}{2}.
      \end{array}         \right.
      \end{eqnarray*}

       Fig. 2 shows numerical results of  Test. 1 for schemes Q-tra1 and Q-tra2. In this test, we have  chosen $cfl=0.5$ and $t=0.5$ on 200 computational cells. The number of cells is large and the value of t is small, so under these conditions, one can see that the performance of both schemes is similar. 
       
  %%%% figura 2%%%%%%%%%%%%%%%%%%%%%%%%%%%%%%%%%%%%%%%
 \vspace{2 cm}
  \begin{figure}[ht]
   \begin{center}
    \begin{picture}(500,150)
    
     \put(-20,-20){ \includegraphics[width=8.2cm]{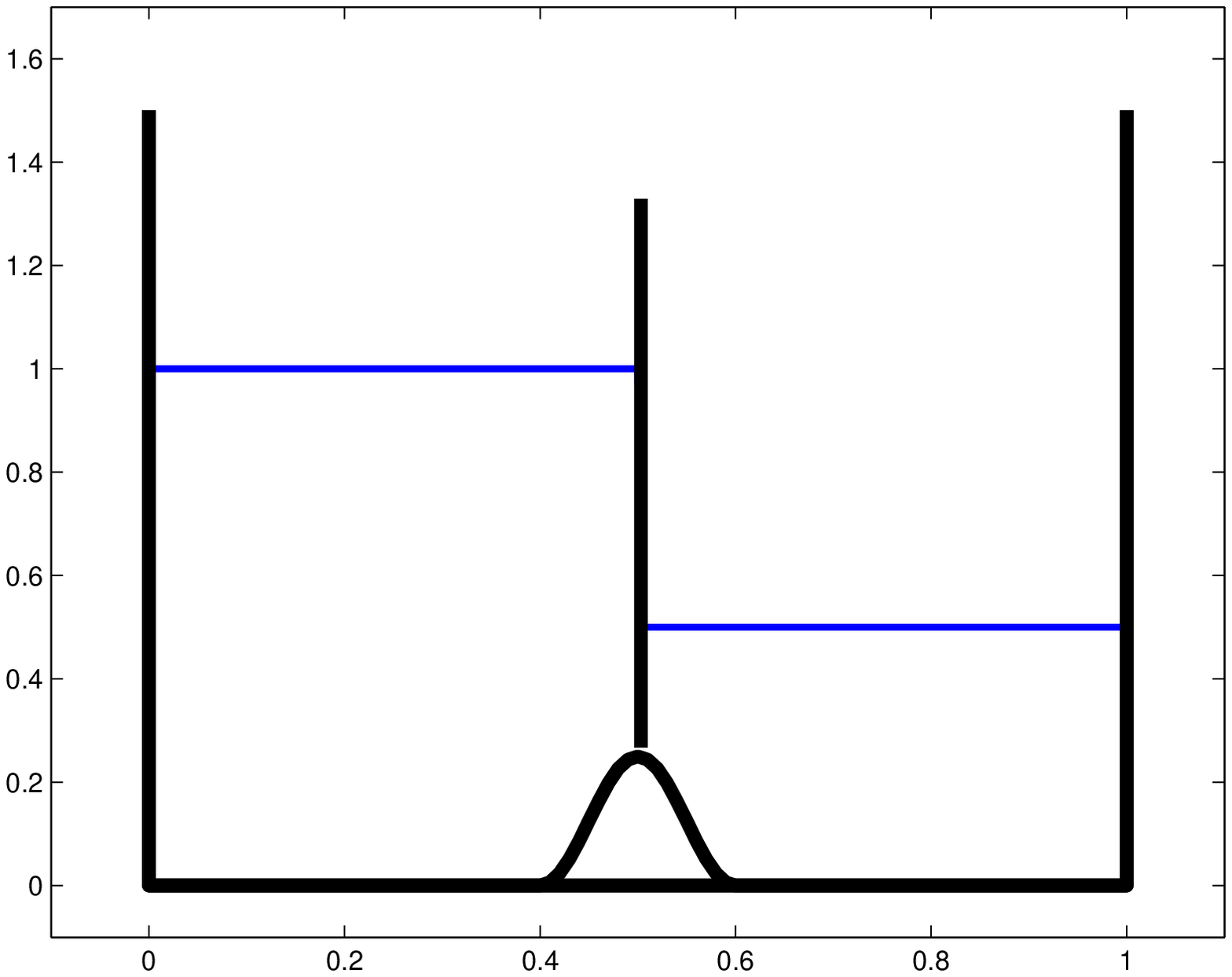} }     
     \put(220,-20){ \includegraphics[width=8.2cm]{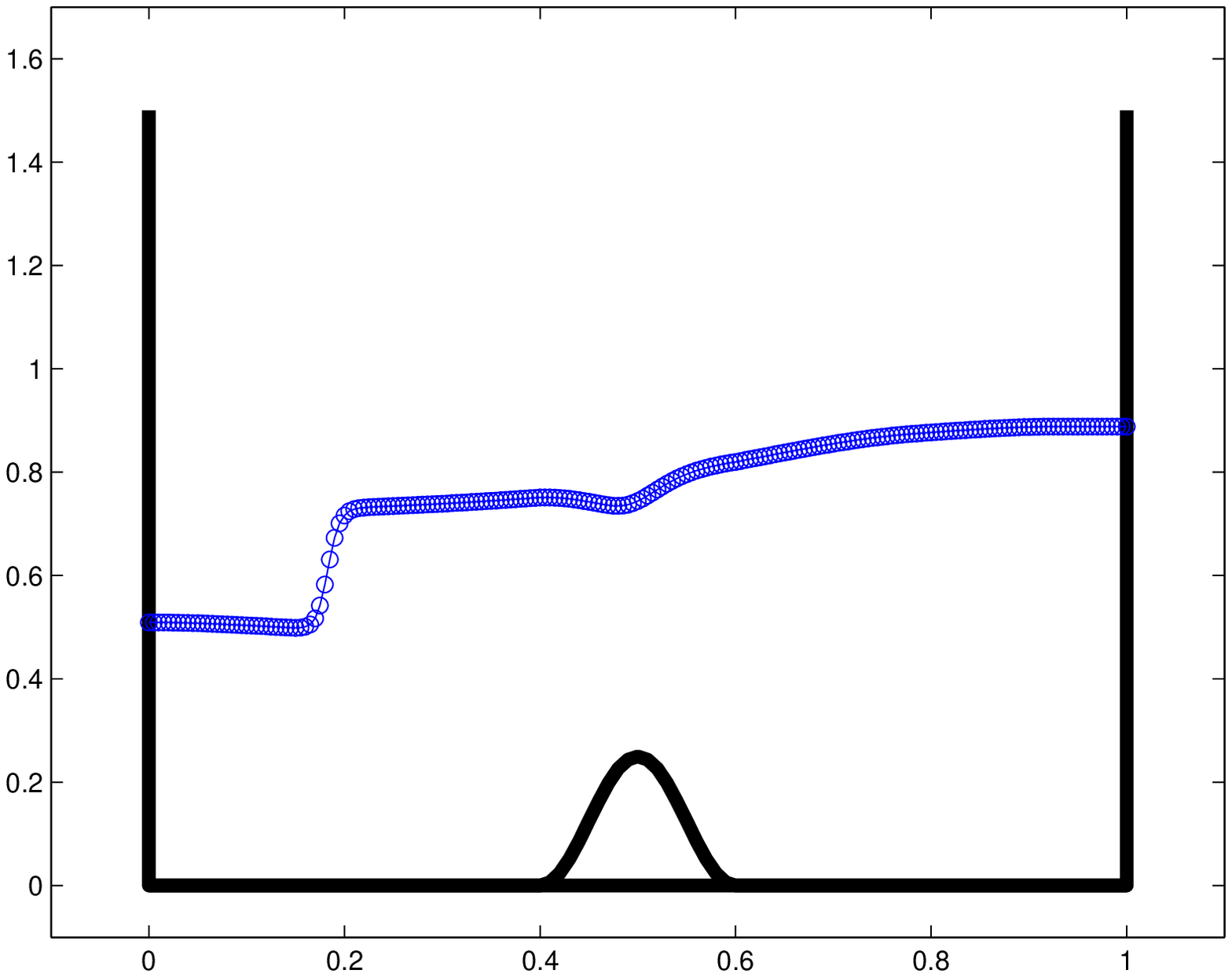} }      
     
      \put(80,150){\footnotesize Initial condition} 
      \put(330,150){\footnotesize $t=0.5$}                   
     
     \put(39,92){\footnotesize free surface}     
     \put(100,117){\footnotesize wall}     
     \put(116,52){\footnotesize free surface}     
     
     \put(320,88){\footnotesize  {\color{blue} -------- } Q-tra1}         
       \put(320,100){\footnotesize {\tiny \color{blue} $\circ$$\circ$$\circ$$\circ$$\circ$$\circ \:\:\:$}  Q-tra2}         
     
        \put(120,17){\footnotesize $b(x)$}
            \end{picture}
    \end{center}
    \caption{Test 1. Initial condition and numerical result for $t=0.5 .$}
    \end{figure}
%%%%%%%%%%%%%%%%%%%%%%%%%%%%%%%%%%%%%%%%%%%%%    

\subsection*{Test 2: Stationary problem with smooth bottom}
\par

We consider also  Eq. (\ref{ec4}) on the domain 
  $(x, t) \in [0,1]  \times \R^+.$  The bottom function $b(x)$ is the same of test 1  and   the initial data
are
 
  \begin{eqnarray*} 
      W(x,0) =             \left( \begin{array}{c}
     1- b(x) \\[0.15 cm]
     0  \end{array} \right) , \: 0 \leq x \leq 1.
      \end{eqnarray*}

This problem has a stationary solution $W(x,t) = W(x,0), \: 0 \leq x \leq 1 \:\: $ and $\: \: t>0.$

In this test, we have chosen $cfl=0.5$ and $t=0.25$ on 50 computational cells, it shows  spurious oscillations for Q-tra1 (see Fig. 3). However,
        Q-tra2 scheme has a good behavior.

        %%%% figura 3%%%%%%%%%%%%%%%%%%%%%%%%%%%%%%%%%%%%%%%
 \vspace{2 cm}
  \begin{figure}[ht]
   \begin{center}
       \vspace*{2 cm}

    \begin{picture}(500,150)
    
     \put(80,-20){ \includegraphics[width=11cm]{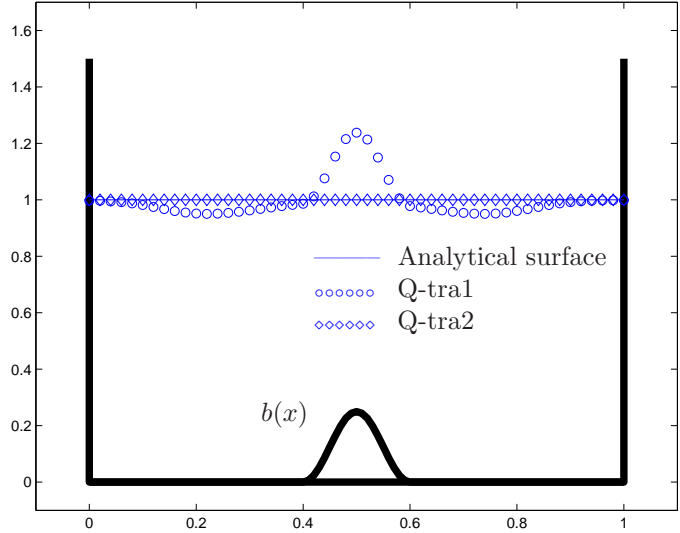} }        
   
         \put(230,99){\footnotesize {\color{blue} -------- } Analytical surface}         
       \put(230,87){\footnotesize {\tiny \color{blue} $\circ$$\circ$$\circ$$\circ$$\circ$$\circ \:\:\:$}  Q-tra1}         
          \put(230,75){\footnotesize {\tiny \color{blue} 
                $\diamond$$\diamond$$\diamond$$\diamond$$\diamond$$\diamond\:\:\: $}  Q-tra2}

        \put(210,40){\footnotesize $b(x)$}
            \end{picture}
    \end{center}
    \caption{Test 2. Numerical result for $t=0.25 .$}
    \end{figure}
%%%%%%%%%%%%%%%%%%%%%%%%%%%%%%%%%%%%%%%%%%%%%    

\subsection*{Test 3: Tidal wave flow with irregular topography}
\par

We consider  Eq. (\ref{ec4}) on the domain 
  $(x, t) \in [0,1500]  \times [0,10800],$ where length is in meters and time in seconds.  The bottom function $b(x)$ is shown in Fig. 4,  and   the initial data
are
 
  \begin{eqnarray*} 
      W(x,0) =             \left( \begin{array}{c}
     H(x) \\[0.15 cm]
     0  \end{array} \right) , \:\: \:  H(x)=H(0)-b(x) , \: \:\:  0 \leq x \leq 1500   \: \:\: 
          \mbox{and}   \: \:\:  H(0)=16.
      \end{eqnarray*}

Moreover, we have the following boundary conditions

     \begin{eqnarray*} 
                \left. 
                 \begin{array}{l}
     h(0,t)= H(0) + \varphi (t), \\[0.35 cm]
     \varphi (t)= \displaystyle 4 +4 \: \sin \left( \pi \left( \frac{4 \: t}{86400}- \frac{1}{2}\right) \right), \\[0.35 cm]
     q(1500, t)= 0 , \end{array} \right\}  \: 0 \leq t \leq 10800.
      \end{eqnarray*}

Function $   \varphi (t) $  simulates a tidal wave of $4 m$ amplitud.

For this problem, we can obtain an asymptotic analytical solution (see \cite{ber2}) given by 
$h(x,t)=H(x)+\varphi (t)$ and it satisfies $h(x,10800)=20 , \: 0 \leq x \leq 1500.$

In order to check the behaviour  of Q-tra2 scheme, we have chosen $cfl=0.9$ on 100 computational cells. 
 Test. 3 is an strong test, it uses a long time $t=10800$ seconds and a complex topography, however 
 an excellent performance of this scheme is demonstrated (see Fig. 4).
 
 %%%% figura 4%%%%%%%%%%%%%%%%%%%%%%%%%%%%%%%%%%%%%%%

  \begin{figure}[ht]
   \begin{center}
    \vspace*{2 cm}
    \begin{picture}(500,180)
    
     \put(80,-20){ \includegraphics[width=11cm]{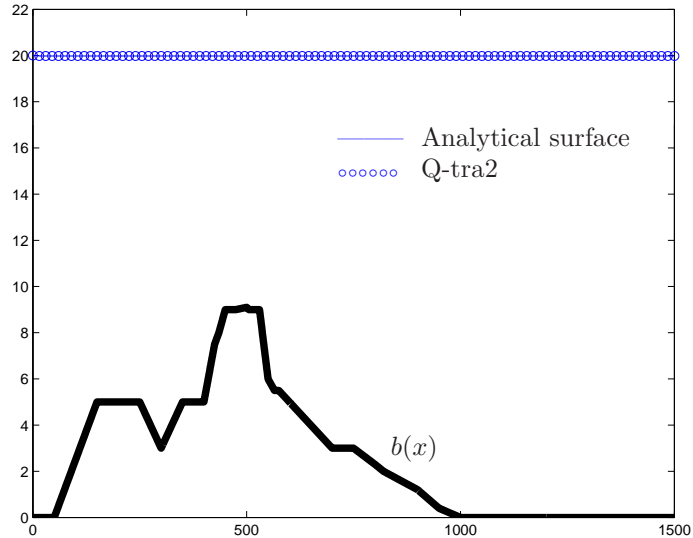} }        
   
         \put(240,135){\footnotesize {\tiny \color{blue} $\circ$$\circ$$\circ$$\circ$$\circ$$\circ \:\:\:$}  Q-tra2}          
          \put(240,147){\footnotesize {\color{blue} -------- } Analytical surface}                  
               
        \put(260,30){\footnotesize $b(x)$}
            \end{picture}
    \end{center}
    \caption{Test 3. Tidal wave flow over an irregular topography using Q-tra2.}
    \end{figure}
%%%%%%%%%%%%%%%%%%%%%%%%%%%%%%%%%%%%%%%%%%%%%    

\subsection*{Test 4: Tidal wave on the shoreline with friction effects}
\par

In this test we take in account the bottom friction.
 Manning law is used to model friction between  fluid and bottom, so a new element appears in the source term. As a consequence, the source term is given by
\begin{eqnarray*} 
\hspace{-0.75 cm}  G(x,W)=   \left( \begin{array}{c}
   \hspace{-0.25 cm}   0 \\[0.15 cm] \displaystyle
    \hspace{-0.25  cm}  - g h b'(x) - g q M^2 \left| \frac{q}{h}  \right| h^{- \frac{4}{3}} 
     \end{array} \right),
\end{eqnarray*}
      where  $M$ is the Manning coefficient.
      
      \par
      
      We consider the domain 
  $(x, t) \in [0,6]  \times \R^+$ and the bottom function    
  \begin{eqnarray*} 
  b(x)= 
      \left\{  \begin{array}{l}
   \displaystyle   0.00125 \, x + 0.0125, \:  0 \leq x \leq 3, \\[0.25 cm]
    0.162 (x-3)+ 0.01625, \:  3 \leq x \leq 6.
     \end{array} \right. 
            \end{eqnarray*}
            
         \vspace{0.75 cm}            
       
        Here, length is in meters, time in seconds and the reference level is located at 40 cm.  
        The tidal wave is simulated by introducing a discharge $q=0.8 m^3/s$ in the boundary $x=0$ during the first $0.2$ seconds, from that time a vertical wall condition is introduced. 
 In the boundary $x=6$ is also introduced a vertical wall condition to preserve the mass conservation in the domain.  
            
\vspace{0.75 cm}
   
   On the other hand, this test involves a wet/dry front, which complicates the stability of the scheme. Specific treatment is required to remove the spurious oscillations, the reader can find all detailed in \cite{bru} and \cite{cas}.  For brevity, we summarize here two main actions of this treatment:
       \begin{description}
        \item[]  \hspace{0.75 cm} 1.- Redefinition of the discretized bottom function to avoid the appearance of 
                     spurious pressure forces.   
     \begin{description}
       \item[] \hspace{1.5 cm} $\bullet$ If $I_j$ is dry, $I_{j-1}$ wet and $h_{j-1}+b_{j-1} < h_j+b_j,$ 
       then $b_j=b_{j-1}+h_{j-1}.$  
         \item[] \hspace{1.5 cm} $\bullet$ If $I_j$ is wet, $I_{j-1}$ dry and $h_{j-1}+b_{j-1} > h_j+b_j,$ 
       then $b_j=b_{j-1}-h_{j}.$  
               \end{description}
        \item[]  \hspace{0.75 cm} 2.-  Simulate the fact that the discharge is zero when the fluid across a wet/dry front.
     \begin{description}
       \item[] \hspace{1.5 cm} $\bullet$ If $I_j$ is dry, then $q^{n+1}_j=0.$  
       \item[] \hspace{1.5 cm} $\bullet$ If $I_j$ is wet, estimate of ${q}^{n+1}_j<0$ and $I_{j-1}$ is dry, then $q^{n+1}_j=0.$         
        \item[] \hspace{1.5 cm} $\bullet$ If $I_j$ is wet, estimate of ${q}^{n+1}_j>0$ and $I_{j+1}$ is dry, then $q^{n+1}_j=0.$ 
      \end{description}
        \end{description}

   \vspace{0.75 cm}
      
      In addition, friction causes spurious oscillations in  $q^{n}_j.$  However splitting schemes can avoid these effects by using a semi-implicit discretization.  Below, we summarize the details.

   \vspace{0.75 cm}
       
       In order to compute $W^{n+1}_{Lj},$ we consider Eq. (\ref{ec21}) with friction
       
        \begin{eqnarray} \label{ec30}
     \left\{ \begin{array}{l}
  \displaystyle \frac{d}{dt} \left( \begin{array}{l}
     h_j \\[0.25 cm]
     q_j \end{array} \right) =    \left( \begin{array}{c}
     \displaystyle 
      -  d_{12}^1 g h_j b'_j -  d_{12}^1 g  q_j M^2 \left| \frac{q_j}{h_j} \right| \left( h_j \right)^{- \frac{4}{3}} \\[0.8 cm]
        \displaystyle 
               -  d_{12}^1 g h_j b'_j -  d_{12}^1 g  q_j M^2 \left| \frac{q_j}{h_j} \right| \left( h_j \right)^{- \frac{4}{3}}
        \end{array}  \right); \: \: \: \: \:  t\in [t^n, t^{n+1}] , \\[2 cm]
      \left( \begin{array}{l}
     h(x_j,t^n) \\[0.15 cm]
     q(x_j, t^n) \end{array} \right) =       
      \left( \begin{array}{l}
     \hat{h}(x_j,t^{n+1}) \\[0.25 cm]
     \hat{q}(x_j, t^{n+1}) \end{array} \right)      .
     \end{array} \right.
      \end{eqnarray} 
       
     \vspace{0.75 cm}
       Solving, we have 
       \begin{equation}\label{ec31}
                      h_{L_j}=   
   \hat{h}_j - \displaystyle  d^1_{12} \Delta t g \frac{h_{j}+h_{j-1}}{2} b'_j
         -  d^1_{12} \Delta t g \left( \frac{h_{j}+h_{j-1}}{2} \right)^{-\frac{4}{3}} 
      \frac{q_{j}+q_{j-1}}{2} M^2 \left| \frac{q_{j}+q_{j-1}}{h_{j}+h_{j-1}}  \right|
          \end{equation}
and         
          \begin{equation}\label{ec32}
                     q_{L_j}=   
   \hat{q}_j - \displaystyle  d^1_{22} \Delta t g \frac{h_{j}+h_{j-1}}{2} b'_j
         -  d^1_{22} \Delta t g \left( \frac{h_{j}+h_{j-1}}{2} \right)^{-\frac{4}{3}} 
      q_{L_j} M^2 \left| \frac{q_{j}+q_{j-1}}{h_{j}+h_{j-1}}  \right|.
          \end{equation}
        
      \vspace{0.75 cm} 
       
       Now, from Eq. (\ref{ec32}) we obtain that
       
       \begin{equation}\label{ec33}
                     q_{L_j}=    \displaystyle
                      \frac{  \hat{q}_j - \displaystyle  d^1_{22} \Delta t g \frac{h_{j}+h_{j-1}}{2} b'_j }
                     { 1  +  \displaystyle d^1_{22} \Delta t g \left( \frac{h_{j}+h_{j-1}}{2} \right)^{-\frac{4}{3}} 
 M^2 \left| \frac{q_{j}+q_{j-1}}{h_{j}+h_{j-1}}  \right|}.
          \end{equation}

     \vspace{0.75 cm}
             
 On the other hand,  in order to compute $W^{n+1}_{Rj},$ we act analogously and obtain
 
  \begin{equation}\label{ec34}
                      h_{R_j}=   
   \hat{h}_j - \displaystyle  d^2_{12} \Delta t g \frac{h_{j+1}+h_{j}}{2} b'_{j+1}
         -  d^2_{12} \Delta t g \left( \frac{h_{j+1}+h_{j}}{2} \right)^{-\frac{4}{3}} 
      \frac{q_{j+1}+q_{j}}{2} M^2 \left| \frac{q_{j+1}+q_{j}}{h_{j+1}+h_{j}}  \right|
          \end{equation}
and         
          \begin{equation}\label{ec35}
                     q_{R_j}=   
             \displaystyle
                      \frac{  \hat{q}_j - \displaystyle  d^2_{22} \Delta t g \frac{h_{j+1}+h_{j}}{2} b'_{j+1} }
                     { 1  +  \displaystyle d^2_{22} \Delta t g \left( \frac{h_{j+1}+h_{j}}{2} \right)^{-\frac{4}{3}} 
 M^2 \left| \frac{q_{j+1}+q_{j}}{h_{j+1}+h_{j}}  \right|}.
          \end{equation}
 
\vspace{0.75 cm}

               \begin{remark} 
   We denote the time splitting scheme given by 
   (\ref{ec14})-(\ref{ec31})-(\ref{ec33})-(\ref{ec34}) and (\ref{ec35}) as {\bf Q-tra3}.
  \end{remark}
                   
      \vspace{0.75 cm}
      
            Test. 4 results using Q-tra3 scheme are shown in Fig. 5. This test has all the ingredients of a real problem in shallow water: irregular topography, friction effects between bottom and water, and advancing front wet/dry. In Fig. 5, we can see                      
           the progress of the tidal wave for different time values: for  $t = 1$ the wave goes to head, in $t=2$ the wave   
     reaches the shoreline,  in $t=3$ the wave hits the wall and in $t=4, $ it returns toward $x = 0.$ In $t=5,$     
     we can see how the wave returns back with negative velocity and it crashes with the wave from $x=0$ with
     positive velocity. However, in front of this critical situation, the behavior of the scheme is stable.

     In this test,
     we have used a Manning coefficient  $M=0.015 \:$ and $\: cfl=0.5$ on 250 computational cells.
        It is appropriate to note that if you do not use a semi-implicit scheme, we obtain spurious oscillations and
         overflow computations.

        %%%% figura 5%%%%%%%%%%%%%%%%%%%%%%%%%%%%%%%%%%%%%%%
 \vspace{2 cm}
  \begin{figure}[ht]
   \begin{center}
    \begin{picture}(500,350)
    
         \put(-20,170){ \includegraphics[width=8.2cm]{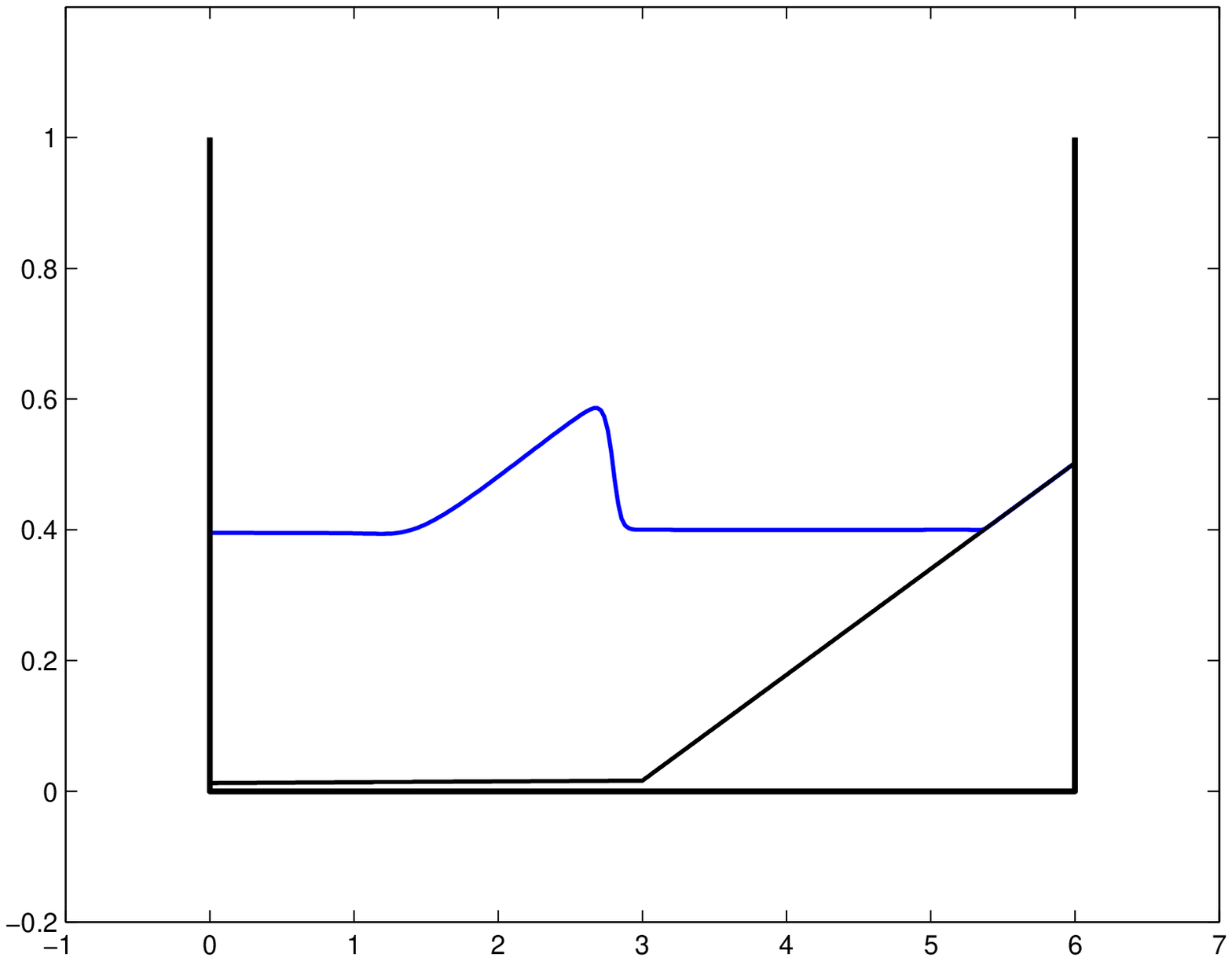} }     
         \put(220,170){ \includegraphics[width=8.2cm]{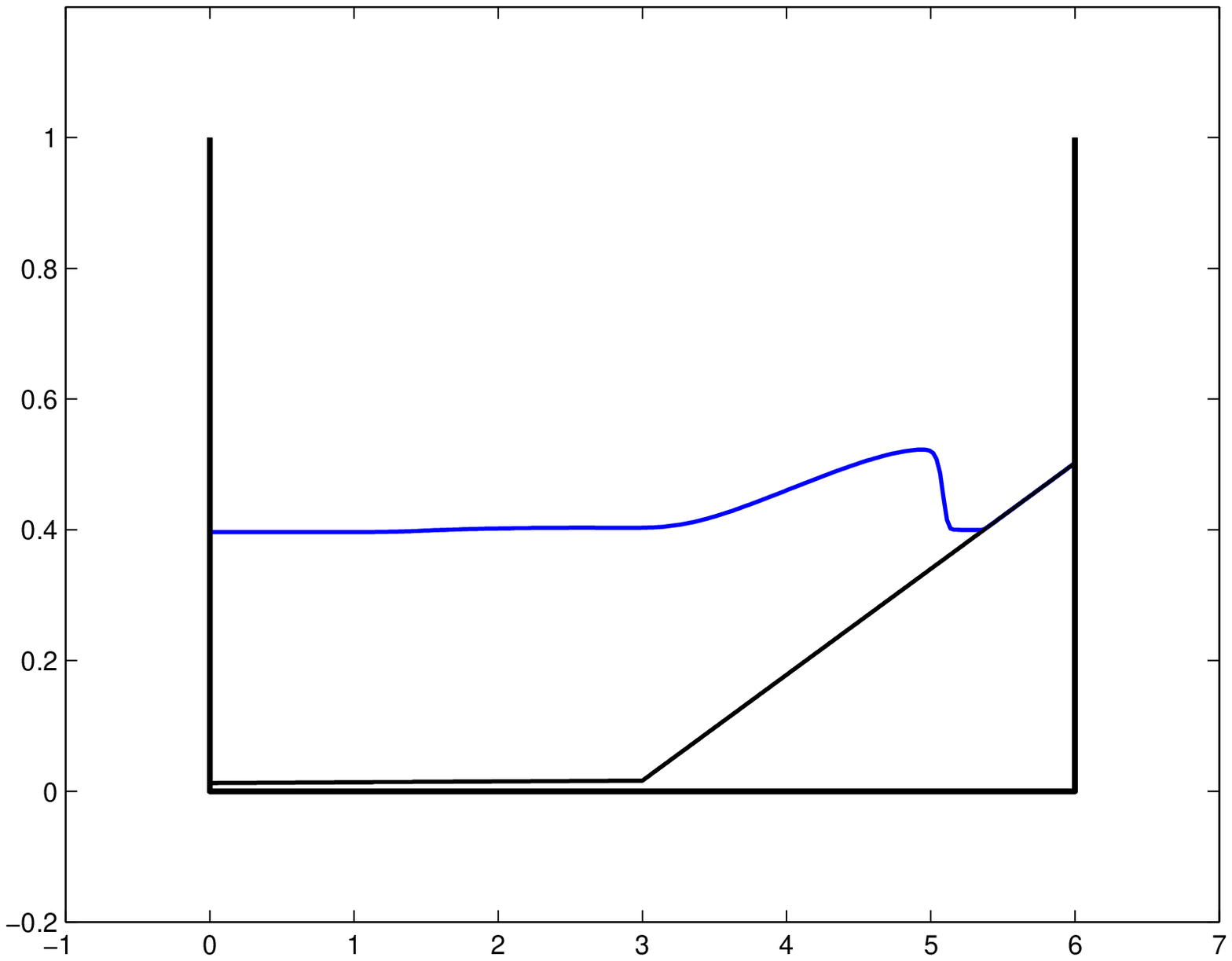} }      
         
         \put(100,340){\footnotesize $t=1$} 
        \put(330,340){\footnotesize $t=2$}  
        \put(60,280){\footnotesize free surface}     
     \put(160,320){\footnotesize wall}     
       \put(140,227){\footnotesize $b(x)$}
         \put(310,265){\footnotesize free surface}     
     \put(400,320){\footnotesize wall}     
       \put(380,227){\footnotesize $b(x)$}

      \put(-20,-20){ \includegraphics[width=8.2cm]{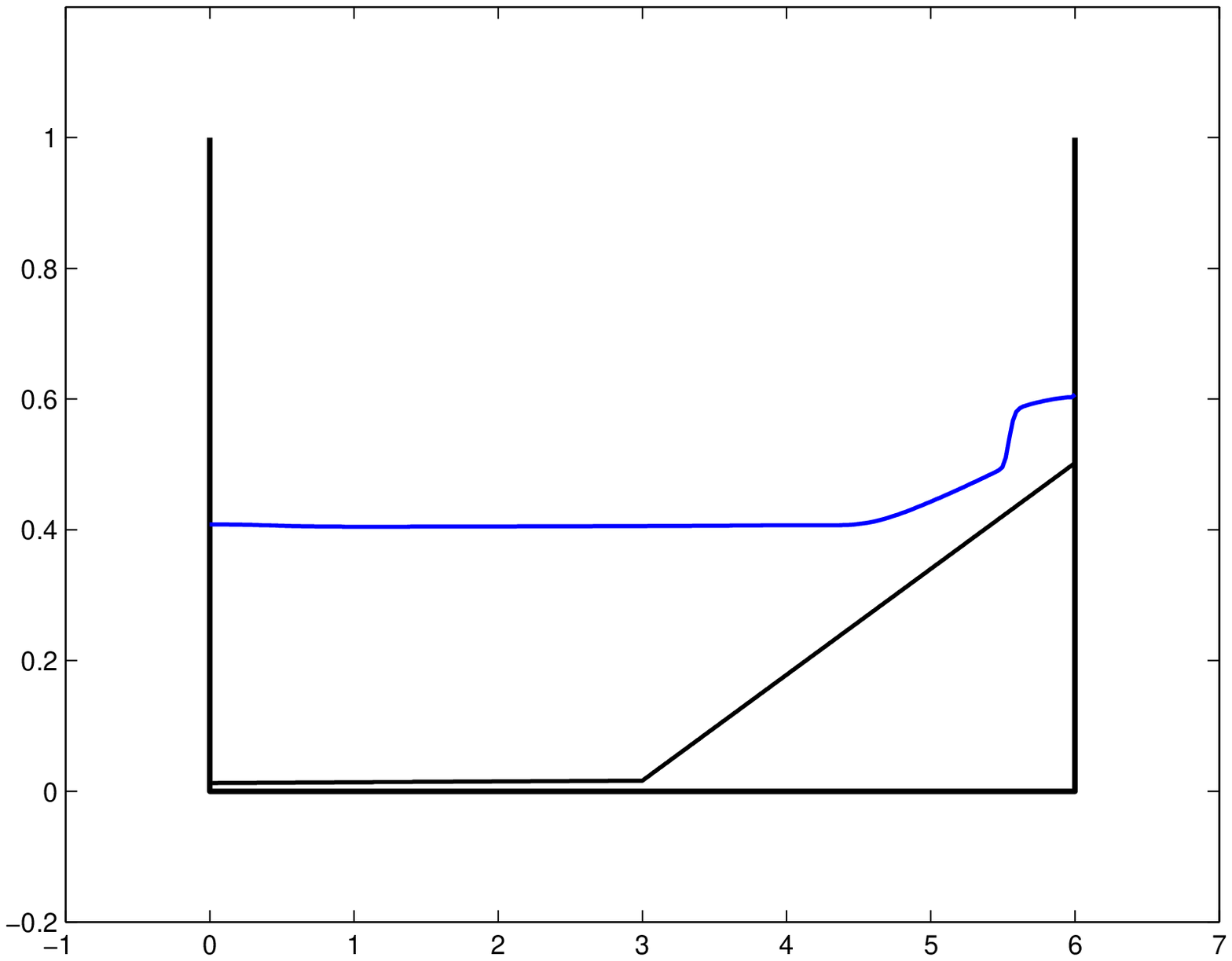} }     
     \put(220,-20){ \includegraphics[width=8.2cm]{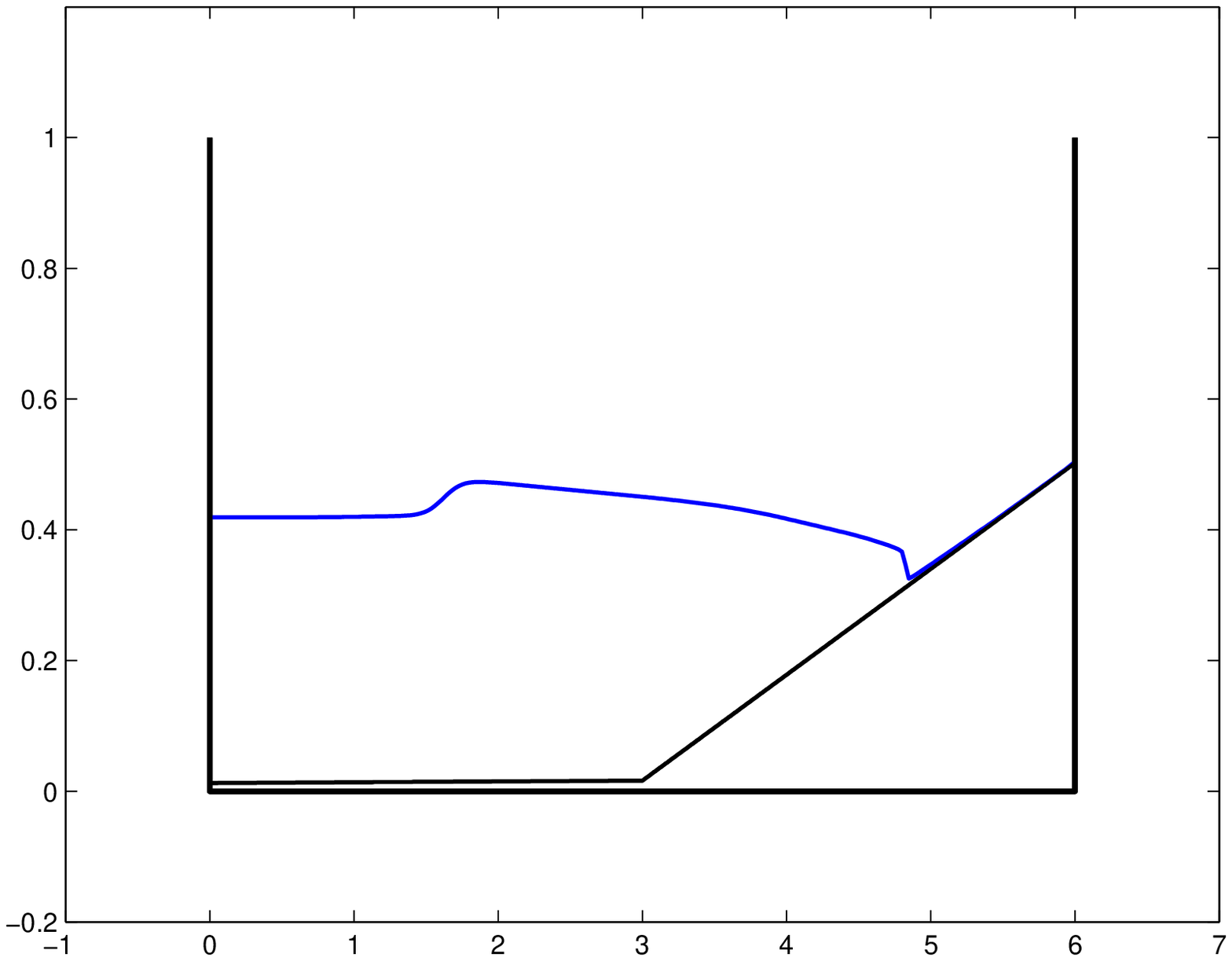} }      
     
      \put(100,150){\footnotesize $t=3$} 
      \put(330,150){\footnotesize $t=5$}                   
     \put(50,72){\footnotesize free surface}     
     \put(160,130){\footnotesize wall}     
       \put(125,27){\footnotesize $b(x)$}
         \put(310,75){\footnotesize free surface}     
     \put(400,130){\footnotesize wall}     
       \put(370,27){\footnotesize $b(x)$}
       
            \end{picture}
                \end{center}
    \caption{Test 4. Tidal wave evolution for different times using Q-tra3.}
    \vspace{2 cm}
        \end{figure}
%%%%%%%%%%%%%%%%%%%%%%%%%%%%%%%%%%%%%%%%%%%%%    

  \section{Concluding remarks}

Some  upwind time splitting  schemes have been developed for solving 1-D shallow water equations.
In some test problems,  our analysis shows that it is not enough to verify an approximate C-property to obtain good
approximations; we need
that an exact C-property to be verified. Perhaps, for some undemanding tests, it is not necessary, but when the test is demanding, then schemes that only satisfy an approximate C-property can bring to numerical instabilities.

We have shown that time splitting schemes for solving shallow water equations must be well balanced
in the following sense: the scheme used to solve the homogeneous equation Eq. (\ref{ec2}) and the scheme used to integrate the time-dependent ODE Eq. (\ref{ec3}) cannot be whatever. They must be linked so that 
an exact C-property is satisfied. Apparently, Q-tra2 is computationally very expensive. But this is not true, since carried out calculations to obtain the numerical approximation of Eq. (\ref{ec2}) are used to obtain $D_L$ and $D_R$ matrices.

Moreover, time splitting schemes could give us the possibility to choose a semi-implicit or implicit solver in order to
obtain stability when Eq. (\ref{ec3})  is integrated. In Test. 4, we have used a semi-implicit scheme that allowed us to
obtain stability when shallow water model takes into account friction effects between fluid and bottom.

In summary, we have presented three time splitting schemes with different properties, the second of which is the best suited to treat severe test on shallow water equations. The semi-implicit Q-tra3 scheme shows the versatility that can provide splitting techniques for solving challenges of stiffness that present some partial differential equations.

 \vspace*{1.5 cm}
 
 \section*{Acknowledgements}

  This work was supported
  by the plan for promotion of research at the University Jaume I, Project P1$\cdot$1B2009-55.
   The author is pleased to acknowledge Professor Fernando Casas for helpful suggestions.

\end{document}